\documentclass[11pt,fleqn]{cmes-az}
\cmes{1}{1}{2009}{0}

\graphicspath{{./image/}}
\runningtitle{A developed new algorithm for evaluating Adomian polynomials}

\title{A developed new algorithm for evaluating Adomian polynomials}
\author{
 M. Azreg-A\"{\i}nou\thanks{Ba\c{s}kent University, Department of Mathematics, Ba\u{g}l\i ca Campus, Ankara, Turkey.}}

\begin{document}
\maketitle

\abstract{
Adomian polynomials (AP's) are expressed in terms of new objects called reduced polynomials (RP's). These new objects, which carry two subscripts, are independent of the form of the nonlinear operator. Apart from the well-known two properties of AP's, curiously enough no further properties are discussed in the literature. We derive and discuss in full detail the properties of the RP's and AP's. We focus on the case where the nonlinear operator depends on one variable and construct the most general analytical expressions of the RP's for small values of the difference of their subscripts. It is shown that each RP depends on a number of functions equal to the difference of its subscripts plus one. These new properties lead to implement a dramatically simple and compact Mathematica program for the derivation of individual RP's and AP's in their general forms and provide useful hints for elegant hand calculations of AP's. Application of the program is considered.
}

\section{Introduction}
Decomposition methods continue to develop and gain ground in applied mathematics and integral methods. They have been reviewed, modified and applied to different fields of science and engineering. The idea behind a decomposition method, which is used to solve differential, integro-differential, algebraic equations and their combinations, is to obtain by easily handled successive iterations or recursions approximate solutions within a predefined accuracy. From this point of view, the most applied methods are a variety of domain decomposition methods, which solve boundary value problems~[\cite{Davies}; \cite{Han}; \cite{Huang}; \cite{Patr}; \cite{Takei}; \cite{Vod}], and Adomian decomposition method and its modifications, which tackle both initial and boundary values problems~[\cite{ad}; \cite{ad2}; \cite{coch}; \cite{taw}; \cite{hoss}; \cite{Lai}; \cite{waz00}; \cite{waz-say}; \cite{zh}].

Adomian decomposition method and some of its modifications and extensions are mathematical tools providing analytical and rapidly convergent (if not exact) solutions to a variety of problems in nonlinear science. A problem in nonlinear science is modeled by
\begin{equation}\label{1}
    Lu + Ru + Nu = g\,,
\end{equation}
where the unknown (scalar, vector or matrix) function $u(x)$ is subject to some initial and/or boundary conditions. The function $g(x)$ is a source term, ($L+R$) is a linear operator and $N$ includes nonlinear operator terms, if there are any. The original idea due to Adomian is to split the linear operator ($L+R$) into two linear terms, $L$ and $R$, where $L$ represents the highest order derivative and is easily invertible and $R$ groups the remaining lower order derivatives.

Integrating~\eqref{1} leads to a Volterra integral equation~[\cite{coch}]
\begin{equation}\label{2}
    u = f - L^{-1}(Ru) - L^{-1}(Nu)\,,
\end{equation}
where $f$ is the sum of $L^{-1}(g)$ and the terms arising from the application of the initial and/or boundary conditions to $u$. It seems that a part of the method originated in the work of Cochran~[\cite{coch}]. The other new part of the method consists in replacing $u$ by sum of components $u_m$
\begin{equation*}
    u = \sum_{m=0}^\infty u_m
\end{equation*}
leading to
\begin{equation}\label{3}
    \sum_{m=0}^\infty u_m = f - L^{-1}\bigg(R\sum_{m=0}^\infty u_m\bigg) - L^{-1}\bigg(N\sum_{m=0}^\infty u_m\bigg)\,.
\end{equation}
Now, in order to solve~\eqref{1} all one needs is to fix the components $u_m$ using~\eqref{3}. This can be done recursively, however, upon splitting first the action of $N$ on $\sum_{m=0}^\infty u_m$ into sum of terms $\sum_{m=0}^\infty A_m$ in such a way that the first term $A_0$ depends only on $u_0$, which will be the first component of $u$ to be fixed, and $A_1$ depends only on ($u_0,u_1$), ..., and $A_m$ depends only on ($u_0,u_1,\dotsc, u_m$). With that said, the action of $N$ on $u$ reads as
\begin{equation}\label{4}
    N\sum_{m=0}^\infty u_m = \sum_{m=0}^\infty A_m(u_0,u_1,\dotsc, u_m)\,,
\end{equation}
leading to
\begin{equation}\label{5}
    \sum_{m=0}^\infty u_m = f - L^{-1}\bigg(\sum_{m=0}^\infty Ru_m\bigg) - L^{-1}\bigg(\sum_{m=0}^\infty A_m\bigg)\,.
\end{equation}

There is no unequivocal way to fix $u_0$~[\cite{waz-say}], however, a straightforward way to do it is to identify $u_0(x)$ with $f(x)$. Similarly and in a straightforward way the other components are chosen to be
\begin{equation}\label{5}
    u_{m+1} = - L^{-1}(Ru_m) - L^{-1}(A_m),\;m\geq 0,\;\text{with}\;u_0=f\,.
\end{equation}

If the action of $N$ on $u$ is represented by some function $F(u)$, then $A_m(u_0,u_1,\dotsc, u_m)$, for $m>0$, is expressed as a sum of $m$ terms, each term includes a derivative of order $k$ ($1\leq k\leq m$) of $F$ at $u_0$, $F^{(k)}(u_0)$, times a polynomial in $(u_1,u_2,\dotsc, u_m)$ and $A_0=F(u_0)$. The first five so called Adomian ``polynomials", $A_0,\dotsc, A_4$, available in the literature are
\begin{align}
    A_0 &= F(u_0)\,,\nonumber\\
    A_1 &= u_1F^\prime(u_0)\,,\nonumber\\
\label{list1}A_2 &= u_2F^\prime(u_0) + \frac{u_1^2}{2!}F^{\prime\prime}(u_0)\,,\\
    A_3 &= u_3F^\prime(u_0) + u_1u_2F^{\prime\prime}(u_0) + \frac{u_1^3}{3!}F^{(3)}(u_0)\,,\nonumber\\
    A_4 &= u_4F^\prime(u_0) + \bigg(\frac{u_2^2}{2!}+u_1u_3\bigg)F^{\prime\prime}(u_0) + \frac{u_1^2u_2}{2!}F^{(3)}(u_0) + \frac{u_1^4}{4!}F^{(4)}(u_0)\,.\nonumber
\end{align}

Adomian decomposition method (ADM) is known as a reliable mathematical tool for solving algebraic, differential or integral linear and nonlinear equations~[\cite{taw}; \cite{hoss}; \cite{N1}; \cite{Lai}; \cite{N2}; \cite{waz-say}; \cite{zh}] to mention but a few. The method solves nonlinear equations appealing neither to linearization nor to perturbation approaches. ADM provides in many cases the exact solution, if it exists in closed form, or a semianalytical solution showing the main features of the exact one. However, a fastidious task in ADM is the evaluation of AP's which necessitates a big deal of calculations. A first formal formula to evaluate AP's is that introduced by Adomian himself~[\cite{ad}] followed by a series of some other alternative techniques developed in~[\cite{W1}; \cite{W2}] and some other papers aiming to reduce the size of calculations, however, they have fallen short of their aim.  Based on the property that in each $A_m$ the sum of subscripts of the components $(u_0,u_1,\dotsc, u_m)$ is equal to $m$, a new algorithm has been developed to evaluate AP's without using Adomian formula~[\cite{waz00}]. Besides the above mentioned property, AP's have some other interesting properties, never discussed in the literature, which will be discussed and derived in the next section. The aim of this paper is to use these new properties as well as the above mentioned one to derive new reliable techniques for hand and computer evaluations of AP's. From the latter point of view a very simple and compact Mathematica program to compute any polynomial $A_m$ in its most general form as a function of $F^{(k)}(u_0)$, $1\leq k\leq m$, and of $(u_1,u_2,\dotsc, u_m)$ is implemented. When the function $F$ is given, the program derived below allows immediate evaluation of any individual polynomial $A_m$ without evaluating the polynomials with subscripts smaller than $m$ and without using Adomian formula. Furthermore, the use of the new properties allows for elegant hand calculations of AP's.

In section 2 we introduce the reduced polynomials $Z_{m,k}$ and derive and discuss in detail their properties as well as the properties of $A_m$. In section 3 we apply our results to derive general analytical expressions of $A_m$ for small values of $m-k$, which are valid for any $F$. In section 4 we discuss in detail the implementation of the Mathematica program and apply it to solve the problem of the pendulum in section 5. We conclude in section 6.

\section{Properties of Adomian polynomials}
The polynomials $A_m$, $m>0$, possess the following properties
\begin{itemize}
  \item \textit{Property 1.} $A_m$ depends by construction only on the vector $(u_1,u_2,\dotsc,$ $u_m)$ and does not depend on $u_n$ with $n>m$;
  \item \textit{Property 2.} $A_m$ is the sum of $m$ terms of the form
  \begin{equation}\label{Am}
    A_m = \sum_{k=1}^m Z_{m,k}(u_1,u_2,\dotsc, u_m)F^{(k)}(u_0)
  \end{equation}
  where $Z_{m,k}$ are called here the reduced polynomials and depend on $(u_1,u_2,\dotsc,$ $u_m)$. We will see below that the properties of $A_m$ are encoded in $Z_{m,k}$;
  \item \textit{Property 3.} In each reduced polynomial $Z_{m,k}$, the components of the vector $(u_1,u_2,\dotsc,$ $u_m)$ appear in such a way that the sum of their subscripts is equal to $m$~[\cite{waz00}]. A consequence of this property is that each reduced polynomial $Z_{m,k}$ depends only on few components of the vector $(u_1,u_2,\dotsc, u_m)$. This is obvious since the sum of the subscripts of the components of the vector $(u_1,u_2,\dotsc, u_m)$ exceeds $m$;
  \item \textit{Property 4.} Another quit important property, curiously enough it has never been discussed in the literature to our knowledge, is that in each $Z_{m,k}$ the sum of the powers of the components $u_n$, $1\leq n\leq m$, is equal to $k$ (the order of the derivative in~\eqref{Am}). Hence, both subscripts of $Z_{m,k}$ have a direct meaning: one is the sum of the subscripts and the other is the sum of the powers of the components of the vector ($u_1,u_2,$ $\dotsc, u_m$) which appear in the expression of $Z_{m,k}$.
\end{itemize}
The other properties are discussed below. Using the Taylor expansion for $F(u_0+\epsilon)$ in powers of $\epsilon$,  one obtains
\begin{equation}\label{exp}
    F(u)=F\bigg(u_0+\sum_{i=1}^\infty u_i\bigg)=\underbrace{F(u_0)}_{A_0}+\sum_{k=1}^\infty \frac{F^{(k)}(u_0)}{k!}\,(u_1+u_2+\dotsb)^k\,.
\end{equation}
Now, substituting the multinomial formula
\begin{equation*}
    (u_1+u_2+\dotsb)^k = k!\,\sum_{n_1+n_2+\dotsb \,=\,k} \,\frac{\prod_{j=1} u_j^{n_j}}{\prod_{j=1} n_j!}
\end{equation*}
into Eq.~\eqref{exp} and using the properties 3. and 4. of $Z_{m,k}$ ($m>0$), one obtains
\begin{align}
\label{zmk} & Z_{m,k} = \sum_S \frac{u_1^{n_1}u_2^{n_2}\dotsb u_m^{n_m}}{n_1!n_2!\dotsb n_m!} \nonumber \\
& \text{with}\quad S=\bigg\{(n_1,\dots,n_m) \in \mathbb{N}^m \mid \sum_{i=1}^m n_i = k \;\text{and}\;\sum_{i=1}^m i\,n_i = m\bigg\}\,.
\end{align}
In the above expression of $Z_{m,k}$, if some power $n_i=0$ we have then $u_i^{n_i}=1$ and the corresponding component $u_i$ will not appear in the expression of $Z_{m,k}$. Now, since $Z_{m,k}$ depends only on few components of the vector $(u_1,u_2,\dotsc, u_m)$ as mentioned earlier, some of the powers $n_i$ in~\eqref{zmk} must vanish. The summation in~\eqref{zmk} being restricted to the set $S$, the powers $n_i$ are solutions to the system of diophantine equations
\begin{equation}\label{dio}
    n_1+n_2+n_3+\dots +n_m = k\quad \text{and}\quad n_1+2n_2+3n_3+\dots +mn_m = m\,,
\end{equation}
where we are interested to find solutions in $\mathbb{N}$. Subtracting the first equation in~\eqref{dio} from the second one, one obtains
\begin{equation}\label{dio2}
    n_2+2n_3+3n_4+\dots +(m-k)n_{m-k+1}+\dotsb +(m-1)n_m = m-k\,.
\end{equation}
Since all the coefficients in~\eqref{dio2} and the powers $n_i$ are positive integers or zero, Eq.~\eqref{dio2}, and consequently the system~\eqref{dio}, holds only if all the components $n_i$, $i>m-k+1$, vanish identically for all $m$ and $1\leq k\leq m$. Said otherwise, all possible solutions to the system~\eqref{dio} are such that:
\begin{equation}
    n_{m-k+2}=n_{m-k+3}=n_{m-k+4}=\dotsb =n_m \equiv 0\,,
\end{equation}
meaning that $Z_{m,k}$ does never depend on $(u_{m-k+2},u_{m-k+3},\dotsc, u_m)$. We have thus established the proof that $Z_{m,k}$ depends at most on $(u_1,u_2,\dotsc, u_{m-k+1})$, reducing Eq.~\eqref{zmk} to
\begin{align}
\label{zmk1} & Z_{m,k} = \sum_{S_R} \frac{u_1^{n_1}u_2^{n_2}\dotsb u_{m-k+1}^{n_{m-k+1}}}{n_1!n_2!\dotsb n_{m-k+1}!}\quad \text{with} \nonumber \\
& S_R=\bigg\{(n_1,\dots,n_{m-k+1}) \in \mathbb{N}^{m-k+1} \mid \sum_{i=1}^{m-k+1} n_i = k \;\text{and}\;\sum_{i=1}^{m-k+1} i\,n_i = m\bigg\}.
\end{align}
Said otherwise, the system~\eqref{dio} is equivalent to the following reduced system of diophantine equations
\begin{equation}\label{dio3}
    n_1+n_2+\dots +n_{m-k+1} = k\quad \text{and}\quad n_1+2n_2+\dots +(m-k+1)n_{m-k+1} = m\,,
\end{equation}
where Eq.~\eqref{dio2} becomes
\begin{equation}\label{dio4}
    n_2+2n_3+3n_4+\dots +(m-k)n_{m-k+1} = m-k\,.
\end{equation}

Does $Z_{m,k}$, as provided by~\eqref{zmk1}, depend on the whole vector $(u_1,u_2,\dotsc,$ $u_{m-k+1})$? We will shortly show that this is the case if $1<k\leq m$, however, if $k=1$, $Z_{m,1}$ does not depend on the whole vector $(u_1,u_2,\dotsc,$ $u_m)$; rather, it depends only on $u_m$.

Assume that $1<k\leq m$. It is straightforward to show that $Z_{m,k}$ depends on $u_1$ and linearly on $u_{m-k+1}$ for $1< k<m$. In fact, if we set $n_{m-k+1}=1$ in~\eqref{dio4} then the remaining powers in~\eqref{dio4} must vanish and this leads using the first equation in~\eqref{dio} to $n_1=k-1$. Furthermore, the solution $n_1=k-1$ and $n_{m-k+1}=1$ (all the other powers are zero) is the unique solution to~\eqref{dio3} for which $n_{m-k+1}\neq 0$. Consequently, $Z_{m,k}$ will always have the term
\begin{equation}\label{term}
    \frac{u_1^{k-1}\,u_{m-k+1}}{(k-1)!}
\end{equation}
in its expansion. For $k=m$, the term~\eqref{term} reduces to $u_1^m/m!$. To show that $Z_{m,k}$ depends on $u_2$ and $u_{m-k}$ we first rewrite the numerator in~\eqref{term} as
\begin{equation*}
    u_1^{k-2}u_1u_{m-k+1}\,.
\end{equation*}
Now, rising the subscript of $u_1$ by 1 and lowering the subscript of $u_{m-k+1}$ by 1, in the above expression, while keeping $u_1^{k-2}$ unchanged, leads to
\begin{equation*}
    u_1^{k-2}u_1u_{m-k+1} \rightarrow u_1^{k-2}u_2u_{m-k}\,.
\end{equation*}
Notice that such a transformation does not affect both the sum of the powers and the sum of the subscripts. Hence, $n_1=k-2$, $n_2=1$ and $n_{m-k}=1$ is a solution to~\eqref{dio3} (or~\eqref{dio}) and
\begin{equation}\label{term2}
    \frac{u_1^{k-2}\,u_2\,u_{m-k}}{(k-2)!}
\end{equation}
is a term in the expansion of $Z_{m,k}$. Continuing so by rising one subscript and lowering the other subscript in~\eqref{term2} while keeping $u_1^{k-2}$ unchanged, leads to
\begin{equation}\label{term3}
    \frac{u_1^{k-2}\,u_i\,u_{m-k-(i-2)}}{(k-2)!}\quad (2\leq i\leq i_{\text{max}})
\end{equation}
where $i_{\text{max}}=(m-k+2)/2$ if $(m-k)$ is even or $i_{\text{max}}=(m-k+1)/2$ if $(m-k)$ is odd. Terms of the form~\eqref{term3} all appear in the expansion of $Z_{m,k}$. Now, as $i$ runs from 2 to $i_{\text{max}}$, the components $u_i$ and $u_{m-k-(i-2)}$ span the whole vector $(u_2,u_3,\dotsc,$ $u_{m-k})$ meaning that $Z_{m,k}$ depends on $(u_2,u_3,\dotsc,$ $u_{m-k})$ and consequently, taking into account~\eqref{term}, it depends explicitly on the whole vector $(u_1,u_2,\dotsc,$ $u_{m-k+1})$.

Assume $k=1$. The system of equations~\eqref{dio3} reduces to
\begin{equation}\label{k1}
    n_1+n_2+\dotsb+n_m = 1\;\;\text{and}\;\;n_1+2n_2+\dotsb+mn_m = m\,.
\end{equation}
Subtracting the second equation from the first one $\times$ $m$, one obtains
\begin{equation*}
    (m-1)n_1+(m-2)n_2+\dotsb (m-i)n_i+\dots +2n_{m-2}+n_{m-1}=0\,,
\end{equation*}
where the left-hand side is a sum of positive numbers with positive coefficients. Hence, the unique solution to the previous equation is $n_1=n_2=\dotsb=n_{m-1}=0$. Substituting into~\eqref{k1} leads to $n_m=1$. We have then shown that the unique solution to~\eqref{k1} is $n_m=1$ and the other powers are zero, so $Z_{m,1}$ depends only on $u_m$ and is given by
\begin{equation}\label{zm1}
    Z_{m,1} = u_m\,.
\end{equation}

From the above discussion one can now write a more or less developed expression for $A_m$ ($m>0$)
\begin{align}\label{Amb}
    A_m &= Z_{m,1}(u_m)F^{\prime}(u_0)+\dotsb+Z_{m,k}(u_1,u_2,\dotsc, u_{m-k+1})F^{(k)}(u_0)\nonumber\\
    \quad & \quad +\dotsb +Z_{m,m-3}(u_1,u_2,u_3,u_4)F^{(m-3)}(u_0)+Z_{m,m-2}(u_1,u_2,u_3)F^{(m-2)}(u_0)\nonumber\\
    \quad & \quad +\,Z_{m,m-1}(u_1,u_2)F^{(m-1)}(u_0)+Z_{m,m}(u_1)F^{(m)}(u_0)\,.
\end{align}
If one reads the formula~\eqref{Amb} from the right to the left one sees that $u_1$ appears for the first time in $Z_{m,m}$, that $u_2$ appears for the first time in $Z_{m,m-1}$, that $u_3$ appears for the first time in $Z_{m,m-2}$, ..., $u_{m-1}$ appears for the first time in $Z_{m,2}$ and $u_m$ appears for the first time in $Z_{m,1}$.

One can still derive further properties of $Z_{m,k}$ from~\eqref{dio4} and~\eqref{dio3} (or~\eqref{dio2} and~\eqref{dio}), however, for the analysis that will follow in the next section we will rely only on the properties discussed here. We can now apply the previous results to derive some expressions for $Z_{m,k}$, $m$ being fixed, which will serve as samples for the implementation the Mathematica program in the next section.

\section{Applications and further properties}
According to~\eqref{Amb} (or~\eqref{zmk1}), $Z_{m,m}$ depends only on $u_1$ and to obtain the value of the power $n_1$ we solve~\eqref{dio3} or \eqref{dio}: $n_1=m$. Hence
\begin{equation}\label{m}
    Z_{m,m} = \frac{u_1^m}{m!}\,.
\end{equation}
Similarly, $Z_{m,m-1}$ depends on $u_1$ and $u_2$. Solving~\eqref{dio3} which reduces to:
\begin{equation*}
    n_1+n_2=m-1\quad \text{and}\quad n_1+2n_2=m
\end{equation*}
leads to $n_2=1$ and $n_1=m-2$. Hence
\begin{equation}\label{m-1}
    Z_{m,m-1} = \frac{u_1^{m-2}u_2}{(m-2)!}\,,
\end{equation}
which is of the form~\eqref{term}. $Z_{m,m-2}$ depends on ($u_1,u_2,u_3$). The powers ($n_1,n_2,n_3$) satisfy the system of equations
\begin{equation}\label{m-2}
    n_1+n_2+n_3=m-2\quad \text{and}\quad n_1+2n_2+3n_3=m\,,
\end{equation}
which leads upon subtracting the first equation from the second one to $n_2+2n_3=2$. This last equation has two possible solutions in $\mathbb{N}$: ($n_2=0, n_3=1$) or ($n_2=2, n_3=0$). We obtain the corresponding values of $n_1$ from the first equation in~\eqref{m-2}: $n_1=m-3$ or $n_1=m-4$, respectively. Hence
\begin{align}
    Z_{m,m-2} &= \frac{u_1^{m-3}u_2^0u_3^1}{(m-3)!0!1!}+\frac{u_1^{m-4}u_2^2u_3^0}{(m-4)!2!0!}\nonumber\\
\label{m-2b} \quad &= \frac{u_1^{m-3}u_3}{(m-3)!}+\frac{u_1^{m-4}u_2^2}{(m-4)!2!}\,.
\end{align}
Similarly, one obtains
\begin{align}
    Z_{m,m-3} &= \frac{u_1^{m-4}u_4}{(m-4)!}+\frac{u_1^{m-5}u_2u_3}{(m-5)!}+
    \frac{u_1^{m-6}u_2^3}{(m-6)!3!}\,;\nonumber\\
    Z_{m,m-4} &= \frac{u_1^{m-5}u_5}{(m-5)!}+\frac{u_1^{m-6}u_2u_4}{(m-6)!}+
    \frac{u_1^{m-6}u_3^2}{(m-6)!2!}+\frac{u_1^{m-7}u_2^2u_3}{(m-7)!2!}+
    \frac{u_1^{m-8}u_2^4}{(m-8)!4!}\,;\nonumber\\
    Z_{m,m-5} &= \frac{u_1^{m-6}u_6}{(m-6)!}+\frac{u_1^{m-7}u_3u_4}{(m-7)!}+
    \frac{u_1^{m-7}u_2u_5}{(m-7)!}+\frac{u_1^{m-8}u_2u_3^2}{(m-8)!2!}+
    \frac{u_1^{m-8}u_2^2u_4}{(m-8)!2!}\nonumber\\
\label{m-3}    \quad & \quad +\,\frac{u_1^{m-9}u_2^3u_3}{(m-9)!3!}+\frac{u_1^{m-10}u_2^5}{(m-10)!5!}\,;\\
    Z_{m,1} &= u_m\,,\quad \text{this has been derived in~\eqref{zm1}}\,.\nonumber
\end{align}
We have also obtained the expression of $Z_{m,m-6}$ which is the sum of 11 terms. One can use the above formulas to evaluate the AP's $A_1$ to $A_7$. For instance, if $m=7$ one obtains $Z_{7,2}$ from the expression of $Z_{m,m-5}$ dropping all negative powers of $u_1$ since they are not in $\mathbb{N}$. Hence, $Z_{7,2}=u_1u_6+u_3u_4+u_2u_5$. Similarly, one obtains the other $Z_{7,k}$ ($1\leq k\leq 7$) from the above formulas.

A final conclusion, which can be read from~\eqref{m-1}, \eqref{m-2b} and~\eqref{m-3} (or from~\eqref{zmk1}), is that each reduced polynomial $Z_{m,k}$ is a sum of some number of terms, each term is a product of $k$ factors and each factor represents one of the components of the vector $(u_1,u_2,\dotsc,$ $u_{m-k+1})$. For instance, $Z_{7,2}$ depends on $7-2+1=6$ functions: $(u_1,u_2,\dotsc,$ $u_6)$. $Z_{7,2}$ is a sum of some terms (3 terms), each term is a product of $k=2$ factors and each factor is one of the six functions $(u_1,u_2,\dotsc,$ $u_6)$ such that the sum of the powers is 2 and the sum of the subscripts is 7. This remark helps, in fact, improving hand calculations of $A_m$.

\section{Implementation of the Mathematica program}
The Mathematica \texttt{Reduce} function helps solving diophantine equations. For instance, Eqs.~\eqref{dio3} for $m=4$ and $k=2$ write as
\begin{equation}\label{e1}
    \sum_{i=1}^3 n_i=2\quad \text{and}\quad \sum_{i=1}^3 in_i=4\,.
\end{equation}
Mathematica solves~\eqref{e1} as follows:
\begin{align*}
  & \pmb{\text{Reduce}[\underset{i=1}{\overset{3}{\sum }}n_i==2 \;\&\& \;\underset{i=1}{\overset{3}{\sum }} i n_i==4\;
  \&\& \;n_1 \geq  0 \;\&\& \;n_2 \geq  0\;\&\& \;n_3 \geq  0} \\
  & \pmb{\&\& \;\{n_1,n_2,n_3\} \in  \text{Integers}]} \\
  & (n_1==0\;\&\&\;n_2==2\;\&\&\;n_3==0)\| (n_1==1\;\&\&\;n_2==0\;\&\&\;n_3==1)
\end{align*}
There are two possible solutions and $Z_{4,2}=(u_2^2/2!)+u_1u_3$ which is the coefficient of $F^{\prime\prime}(u_0)$ in the fifth line of~\eqref{list1}. In the above codes, the conditions $\pmb{n_1 \geq  0}$ $\pmb{\&\&}$ $\pmb{n_2 \geq  0}$ $\pmb{\&\&}$ $\pmb{n_3 \geq  0}$ and $\pmb{\{n_1,n_2,n_3\} \in  \text{Integers}}$ depend on $m$ and are not convenient for large values of $m$. We replace them by the equivalent expressions $\pmb{\underset{i=1}{\overset{3}{\sum }}\text{Abs}[n_i]==\underset{i=1}{\overset{3}{\sum }}n_i}$ and $\pmb{\text{Table}[n_i,\{i,3\}]\in \text{Integers}}$, respectively. These new expressions are very important for the mathematica program we will develop in subsection~\ref{sb3} where only $m$ and $F$ remain as free parameters. In the following two subsections, we provide two further applications of \texttt{Reduce} to solve~\eqref{dio3} and~\eqref{dio} for $m=4$ and $1\leq k\leq 4$, respectively, and introduce and discuss in detail the codes needed for the Mathematica program.

\subsection{Solving Eqs.~\eqref{dio3} for $m=4$}
In the following Mathematica codes we solve~\eqref{dio3} for $m=4$ and $1\leq k\leq 4$, however, the codes are valid for any value of $m$: changing $m$ leads to new results. In the codes $P$ represents the list of all possible power solutions. It is a list of four objects (elements) where each object is made of one or many sublists. For manipulating lists and sublists see~[\cite{Don}].
\begin{align*}
    & \pmb{m=4; \;P = \text{Table}[0,\{j,m\}];\;k = 0;}\\
    & \pmb{\text{While}[(k=k+1)\leq m,\;P[[k]]= \text{Reduce}[\underset{i=1}{\overset{m-k+1}{\sum }}n_i==k \;\&\&
    \;\underset{i=1}{\overset{m-k+1}{\sum }} i n_i==m}\\
    & \pmb{\&\& \;\underset{i=1}{\overset{m-k+1}{\sum}}\text{Abs}[n_i]==\underset{i=1}{\overset{m-k+1}{\sum }}n_i \;\&\& \;\text{Table}[n_i,\{i,m-k+1\}] \in  \text{Integers}]]}\\
    & \pmb{P}\\
    & \{n_4==1\,\&\&\,n_3==0\,\&\&\,n_2==0\,\&\&\,n_1==0,\\
    & (n_1==0\,\&\&\,n_2==2\,\&\&\,n_3==0)\|(n_1==1\,\&\&\,n_2==0\,\&\&\,n_3==1),\\
    & n_1==2\,\&\&\,n_2==1,\\
    & n_1==4\}
\end{align*}
One sees $P$ as a list of four objects, each object appears in one line. The objects $P[[1]]$, $P[[3]]$ and $P[[4]]$
are each made of one sublist, while the object $P[[2]]$ is made of two sublists separated by $\|$. The object $P[[2]]$ corresponds to $m=4$ and $k=2$, which we know has two solutions. A drawback of solving~\eqref{dio3} is that the different sublists $P[[1]]$, $P[[2,1]]$, $P[[2,2]]$, $P[[3]]$ and $P[[4]]$ have different lengths. For instance, $\text{Length}[P[[1]]]=4$, $\text{Length}[P[[2, 1]]]=3$ and $\text{Length}[P[[3]]]=2$. We will see in the next subsection that solving~\eqref{dio} will always lead to uniform sublists with the same length.

\subsection{Solving Eqs.~\eqref{dio} for $m=4$}
In the following Mathematica codes we solve~\eqref{dio} for $m=4$ and $1\leq k\leq 4$, however, the codes are valid for any value of $m$: changing $m$ leads to new results.
\begin{align*}
    & \pmb{m=4; \;P = \text{Table}[0,\{j,m\}];\;k = 0;}\\
    & \pmb{\text{While}[(k=k+1)\leq m,\;P[[k]]= \text{Reduce}[\underset{i=1}{\overset{m}{\sum }}n_i==k \;\&\&
    \;\underset{i=1}{\overset{m}{\sum }} i n_i==m}\\
    & \pmb{\&\& \;\underset{i=1}{\overset{m}{\sum}}\text{Abs}[n_i]==\underset{i=1}{\overset{m}{\sum }}n_i \;\&\& \;\text{Table}[n_i,\{i,m\}] \in  \text{Integers}]]}\\
    & \pmb{P}\\
    & \{n_4==1\,\&\&\,n_3==0\,\&\&\,n_2==0\,\&\&\,n_1==0,\\
    & (n_4==0\,\&\&\,n_3==0\,\&\&\,n_2==2\,\&\&\,n_1==0)\|\\
    & (n_4==0\,\&\&\,n_3==1\,\&\&\,n_2==0\,\&\&\,n_1==1),\\
    & n_4==0\,\&\&\,n_3==0\,\&\&\,n_2==1\,\&\&\,n_1==2,\\
    & n_4==0\,\&\&\,n_3==0\,\&\&\,n_2==0\,\&\&\,n_1==4\}
\end{align*}
Notice that the sublists $P[[1]]$, $P[[2,1]]$, $P[[2,2]]$, $P[[3]]$ and $P[[4]]$ have the same length 4 ($=m$). In order to be able to extract the values of $Z_{4,k}$ ($1\leq k\leq 4$) and evaluate $A_4$ one needs to generate the corresponding list, $U$, of the functions $u_i$. This is done by the code line
$\pmb{U = P \;\text{/.}\; \{n \to  u\}}$
which replaces $n$ by $u$ in the list of powers $P$. The list $U$ has the same properties as the list $P$.

Now, if the object $U[[k]]$ of $U$ (or the object $P[[k]]$ of $P$) has only one sublist (these are the objects $U[[1]]$, $U[[3]]$ and $U[[4]]$ in our example) then
\begin{equation}\label{sub}
    Z_{m,k}=\frac{\prod_{j=1}^m U[[k,j,1]]^{P[[k,j,2]]}}{\prod_{j=1}^m P[[k,j,2]]!}\,,
\end{equation}
and if the object $U[[k]]$ of $U$ has many sublists (this is the object $U[[2]]$ in our example), the number of sublists being $\text{Length}[U[[k]]]$ leads to
\begin{equation}\label{subm}
    Z_{m,k}=\sum_{i=1}^{\text{Length}[U[[k]]]}\;\frac{\prod_{j=1}^m U[[k,i,j,1]]^{P[[k,i,j,2]]}}{\prod_{j=1}^m P[[k,i,j,2]]!}\,.
\end{equation}

\subsection{The Mathematica program}\label{sb3}
In order to distinguish between the case where $U[[k]]$ has only one sublist and the case where it has many sublists we use the condition $\pmb{U[[k,1,1,1]]==0}$ which holds if $U[[k]]$ has only one sublist with $m>1$. The Mathematica program is readily written for the case $m=7$ which has been chosen for illustration:\\

\noindent\(\pmb{m=7; P = \text{Table}[0,\{j,m\}];}\\
\pmb{\text{If}[m > 1, k = 0;}\\
\pmb{\text{While}[(k=k+1)\leq m,}\\
\pmb{P[[k]]= \text{Reduce}[\underset{i=1}{\overset{m}{\sum }}n_i==k \,\&\&\,\underset{i=1}{\overset{m}{\sum }} i n_i==m\,\&\&\, \underset{i=1}{\overset{m}{\sum
}}\text{Abs}[n_i]==\underset{i=1}{\overset{m}{\sum }}n_i \,\&\& }\\
\pmb{\text{Table}[n_i,\{i,m\}] \in  \text{Integers}]];}\\
\pmb{U = P \,\text{/.}\, \{n \to  u\};}\\
\pmb{A_{m,1} = 0; A_{m,2} = 0;}\\
\pmb{\text{Do}[\text{If}[\text{Length}[U[[k,1,1,1]]] == 0, }\\
\pmb{A_{m,1} = A_{m,1} + \dfrac{\underset{j=1}{\overset{m}{\prod }}U[[k,j,1]]^{P[[k,j,2]]}}{\underset{j=1}{\overset{m}{\prod }}(P[[k,j,2]]!)} \text{Derivative}[k][F][u_0],
}\\
\pmb{A_{m,2} = A_{m,2} + \bigg(\underset{i=1}{\overset{\text{Length}[P[[k]]]}{\sum }}\dfrac{\underset{j=1}{\overset{m}{\prod }}U[[k,i,j,1]]^{P[[k,i,j,2]]}}{\underset{j=1}{\overset{m}{\prod
}}(P[[k,i,j,2]]!)}\bigg) \text{Derivative}[k][F][u_0]], }\\
\pmb{\{k, 1, m\}];}\\
\pmb{A_m = A_{m,1} + A_{m,2}]}\\
\pmb{\text{If}[m \text{==}1 ,A_m = u_1 \,\text{Derivative}[1][F][u_0]]}\\
\pmb{\text{If}[m \text{==}0 ,A_m = F[u_0]]}\)
\begin{align}
  A_7  &= u_7 F'(u_0)+(u_3 u_4+u_2 u_5+u_1 u_6) F''(u_0)\nonumber\\
   \quad & \quad +\bigg(\frac{u_2^2 u_3}{2}+\frac{u_1 u_3^2}{2} + u_1 u_2 u_4+\frac{u_1^2 u_5}{2} \bigg) F^{(3)}(u_0)\nonumber\\
   \quad & \quad +\bigg(\frac{u_1 u_2^3}{6} +\frac{u_1^2 u_2 u_3}{2} +\frac{u_1^3 u_4}{6} \bigg) F^{(4)}(u_0)+\bigg(\frac{u_1^3 u_2^2}{12} +\frac{u_1^4 u_3}{24} \bigg) F^{(5)}(u_0)\nonumber\\
\label{list2} \quad & \quad +\frac{u_1^5 u_2}{120} F^{(6)}(u_0)+\frac{u_1^7}{5040} F^{(7)}(u_0)\,.
\end{align}
The constants that appear in the above expression are products of factorials: $5040=7!$, $12=3!2!$, etc. In case of applications, the functions $F$ and $u_0$ have to be defined earlier.

\section{Application of the program: the pendulum}
A solution for the problem of the pendulum can be written down explicitly using elliptic functions, however, the pendulum is used here for mere illustration of the program since it is known to all workers. In the case where there are no friction and applied forces, the differential equation describing the motion of a pendulum is \begin{equation}\label{pen}
    u_{tt}+b\,\sin u = 0\,,
\end{equation}
where $u(t)$ is the angular displacement and $b$ is a positive ``geometric" constant. The initial conditions are $u(0)=a$ and $u_t(0)=0$, $F(u)\equiv \sin u$ and the operator $L^{-1}$ takes the form
\begin{equation}\label{op}
    L^{-1}[f(t)]=\int_0^t\bigg(\int_0^q [f(p)]\,\text{d}p\bigg)\,\text{d}q\,.
\end{equation}
Applying $L^{-1}$ to both sides of~\eqref{pen} and using the prescribed initial conditions, one obtains
\begin{equation}\label{app}
    u(t)=a-L^{-1}[b\sin(u)]\,.
\end{equation}
The problem of the pendulum is solved upon choosing
\begin{equation}\label{ch}
    u_0=a \quad \text{and}\quad u_{m+1}=-L^{-1}[b\,A_m]\,,\;(m\geq 0)\,,
\end{equation}
where $A_m$ are evaluated using the Mathematica program developed in subsection~\ref{sb3}. We have used ten polynomials $A_m$: from $A_0$ to $A_9$ to evaluate eleven components of $u(t)$: from $u_0(t)$ to $u_{10}(t)$. The results of calculations are shown below as produced by Mathematica\\

\noindent\(\pmb{u = \text{Collect}[\underset{i=0}{\overset{10}{\sum }}u_i,t]}\)

\noindent\(a-\dfrac{b t^2 \text{Sin}[a]}{2} +\dfrac{b^2 t^4 \text{Cos}[a] \text{Sin}[a]}{24} +\\
t^6 \bigg(\dfrac{b^3 \text{Sin}[a]}{360} -\dfrac{b^3 \text{Sin}[3 a]}{720} \bigg)+t^8 \bigg(-\dfrac{b^4 \text{Sin}[2 a]}{5040}+\dfrac{17 b^4 \text{Sin}[4 a]}{161280}\bigg)+\\
t^{10} \bigg(-\dfrac{b^5 \text{Sin}[a]}{45360}+\dfrac{13 b^5 \text{Sin}[3 a]}{604800}-\dfrac{31 b^5 \text{Sin}[5 a]}{3628800}\bigg)+\\
t^{12} \bigg(\dfrac{37 b^6 \text{Sin}[2 a]}{17107200}-\dfrac{37 b^6 \text{Sin}[4 a]}{17107200}+\dfrac{691 b^6 \text{Sin}[6 a]}{958003200}\bigg)+\\
t^{14} \bigg(\dfrac{73 b^7 \text{Sin}[a]}{340540200}-\dfrac{b^7 \text{Sin}[3 a]}{3439800}+\dfrac{9557 b^7 \text{Sin}[5 a]}{43589145600}-\dfrac{5461 b^7 \text{Sin}[7
a]}{87178291200}\bigg)+\\
t^{16} \bigg(-\dfrac{313 b^8 \text{Sin}[2 a]}{12573792000}+\dfrac{37 b^8 \text{Sin}[4 a]}{1067489280}-\dfrac{9683 b^8 \text{Sin}[6 a]}{435891456000}+\dfrac{929569
b^8 \text{Sin}[8 a]}{167382319104000}\bigg)+\\
t^{18} \bigg(-\dfrac{47 b^9 \text{Sin}[a]}{20415732480}+\dfrac{381779 b^9 \text{Sin}[3 a]}{100037089152000}-\\
\dfrac{816337 b^9 \text{Sin}[5 a]}{200074178304000}+\dfrac{1441031 b^9 \text{Sin}[7 a]}{640237370572800}-\dfrac{3202291 b^9 \text{Sin}[9 a]}{6402373705728000}\bigg)+\\
t^{20} \bigg(\dfrac{1884343 b^{10} \text{Sin}[2 a]}{6335682312960000}-\dfrac{540809 b^{10} \text{Sin}[4 a]}{1055947052160000}+\dfrac{31786477 b^{10} \text{Sin}[6
a]}{67580611338240000}-\\
\dfrac{138706613 b^{10} \text{Sin}[8 a]}{608225502044160000}+\dfrac{221930581 b^{10} \text{Sin}[10 a]}{4865804016353280000}\bigg)\,.\)\\

The special case with $a=\pi/2$ reduces to
\begin{equation*}
    u(t)=\frac{\pi }{2}-\frac{b t^2}{2}+\frac{b^3 t^6}{240}-\frac{b^5 t^{10}}{19200}+\frac{11 b^7 t^{14}}{13977600}-\frac{211 b^9 t^{18}}{16293888000}\,.
\end{equation*}

One may wonder whether the expression of $u(t)$ for arbitrary $a$ reduces to the linear case for small values of $a$. A series expansion of $u(t,a)$ about the point $a=0$ and for an arbitrary value of $t$ reads as produced by Mathematica\\

\noindent\(\pmb{\text{Series}[u,\{a,0,3\}]}\)

\noindent\(\bigg(1-\dfrac{b t^2}{2}+\dfrac{b^2 t^4}{24}-\dfrac{b^3 t^6}{720}+\dfrac{b^4 t^8}{40320}-\dfrac{b^5 t^{10}}{3628800}+\dfrac{b^6 t^{12}}{479001600}-\dfrac{b^7
t^{14}}{87178291200}+\\
\dfrac{b^8 t^{16}}{20922789888000}-\dfrac{b^9 t^{18}}{6402373705728000}+\dfrac{b^{10} t^{20}}{2432902008176640000}\bigg) a+\\
\bigg(\dfrac{b t^2}{12}-\dfrac{b^2 t^4}{36}+\dfrac{5 b^3 t^6}{864}-\dfrac{13 b^4 t^8}{15120}+\dfrac{1849 b^5 t^{10}}{21772800}-\dfrac{4153 b^6 t^{12}}{718502400}+\dfrac{149473
b^7 t^{14}}{523069747200}-\\
\dfrac{21019 b^8 t^{16}}{1961511552000}+\dfrac{1100627 b^9 t^{18}}{3492203839488000}-\dfrac{5448101 b^{10} t^{20}}{729870602452992000}\bigg) a^3+O[a]^4\,.\)\\\\
It is straightforward to check that the coefficient of $a$ in the above expression is $\cos(\sqrt{b}\,t)$. Hence, to the first order of approximation the solution is $u(t)=a\,\cos(\sqrt{b}\,t)$ which is the exact solution to the linear problem: $u_{tt}+b\,u = 0$ with the same initial conditions.

\section{Conclusion}
For a nonlinear operator $N$ depending on one variable, it was shown that the reduced polynomial $Z_{m,k}$ depends explicitly on $m-k+1$ functions $(u_1,u_2,\dotsc,$ $u_{m-k+1})$ if $k>1$ and that $Z_{m,1}$ depends only on $u_m$. The AP $A_m$ splits into a sum of $m$ terms where each term is the product of $Z_{m,k}$, which is independent of $N$, with a derivative factor, $F^{(k)}(u_0)$, depending on $N$.

General expressions for $Z_{m,k}$ have been derived for small values of $m-k$, hand calculations have been discussed briefly and a simple and compact Mathematica program has been implemented and applied.

For the purpose of hand calculations and analytical studies, we complete the lists~\eqref{list1} and~\eqref{list2} providing general expressions for $A_m$, which is not available in the literature, by adding the polynomials $A_5$, $A_6$, $A_8$ to $A_{10}$ using our program.\\

\noindent\(\pmb{A_5}\)

\noindent\(u_5 F'[u_0]+(u_2 u_3+u_1 u_4) F''[u_0]+\bigg(\dfrac{u_1 u_2^2}{2} +\dfrac{u_1^2 u_3}{2}\bigg) F^{(3)}[u_0]+\dfrac{u_1^3 u_2}{6} F^{(4)}[u_0]\\+\dfrac{u_1^5}{120}
 F^{(5)}[u_0]\,;\)\\

\noindent\(\pmb{A_6}\)

\noindent\(u_6 F'[u_0]+\bigg(\dfrac{u_3^2}{2}+u_2 u_4+u_1 u_5\bigg) F''[u_0]+\bigg(\dfrac{u_2^3}{6}+u_1 u_2 u_3+\dfrac{u_1^2 u_4}{2}\bigg) F^{(3)}[u_0]\\
+\bigg(\dfrac{u_1^2 u_2^2}{4} +\dfrac{u_1^3 u_3}{6}\bigg) F^{(4)}[u_0]+\dfrac{u_1^4 u_2}{24}  F^{(5)}[u_0]+\dfrac{u_1^6}{720}  F^{(6)}[u_0]\,;\)\\

\noindent\(\pmb{A_8}\)

\noindent\(u_8 F'[u_0]+\bigg(\dfrac{u_4^2}{2}+u_3 u_5+u_2 u_6+u_1 u_7\bigg) F''[u_0]\\
+\bigg(\dfrac{u_2 u_3^2}{2} +\dfrac{u_2^2 u_4}{2} +u_1 u_3 u_4+u_1 u_2 u_5+\dfrac{u_1^2 u_6}{2}\bigg) F^{(3)}[u_0]\\
+\bigg(\dfrac{u_2^4}{24}+\dfrac{u_1 u_2^2 u_3}{2} +\dfrac{u_1^2 u_3^2}{4} +\dfrac{u_1^2 u_2 u_4}{2} +\dfrac{u_1^3 u_5}{6} \bigg) F^{(4)}[u_0]\\
+\bigg(\dfrac{u_1^2 u_2^3}{12} +\dfrac{u_1^3 u_2 u_3}{6} +\dfrac{u_1^4 u_4}{24}\bigg) F^{(5)}[u_0]\\
+\bigg(\dfrac{u_1^4 u_2^2}{48} +\dfrac{u_1^5 u_3}{120}\bigg) F^{(6)}[u_0]+\dfrac{u_1^6 u_2}{720}  F^{(7)}[u_0]+\dfrac{u_1^8 }{40320}F^{(8)}[u_0]\,;\)\\

\noindent\(\pmb{A_9}\)

\noindent\(u_9 F'[u_0]+(u_4 u_5+u_3 u_6+u_2 u_7+u_1 u_8) F''[u_0]\\
+\bigg(\dfrac{u_3^3}{6}+u_2 u_3 u_4+\dfrac{u_1 u_4^2}{2} +\dfrac{u_2^2 u_5}{2} +u_1 u_3 u_5+u_1 u_2 u_6+\dfrac{u_1^2 u_7}{2}\bigg) F^{(3)}[u_0]\\
+\bigg(\dfrac{u_2^3 u_3}{6} +\dfrac{u_1 u_2 u_3^2}{2} +\dfrac{u_1 u_2^2 u_4}{2} +\dfrac{u_1^2 u_3 u_4}{2} +\dfrac{u_1^2 u_2 u_5}{2} +\dfrac{u_1^3 u_6}{6}\bigg) F^{(4)}[u_0]\\
+\bigg(\dfrac{u_1 u_2^4}{24} +\dfrac{u_1^2 u_2^2 u_3}{4} +\dfrac{u_1^3 u_3^2}{12} +\dfrac{u_1^3 u_2 u_4}{6} +\dfrac{u_1^4 u_5}{24}\bigg) F^{(5)}[u_0]\\
+\bigg(\dfrac{u_1^3 u_2^3}{36} +\dfrac{u_1^4 u_2 u_3}{24} +\dfrac{u_1^5 u_4}{120}\bigg) F^{(6)}[u_0]\\
+\bigg(\dfrac{u_1^5 u_2^2}{240} +\dfrac{u_1^6 u_3}{720}\bigg) F^{(7)}[u_0]+\dfrac{u_1^7 u_2 }{5040}F^{(8)}[u_0]+\dfrac{u_1^9 }{362880}F^{(9)}[u_0]\,;\)\\

\noindent\(\pmb{A_{10}}\)

\noindent\(u_{10} F'[u_0]+\bigg(\dfrac{u_5^2}{2}+u_4 u_6+u_3 u_7+u_2 u_8+u_1 u_9\bigg) F''[u_0]\\
+\bigg(\dfrac{u_3^2 u_4}{2}+\dfrac{u_2 u_4^2}{2}+u_2 u_3 u_5+u_1 u_4 u_5+\dfrac{u_2^2 u_6}{2}+u_1 u_3 u_6+u_1 u_2 u_7+\dfrac{u_1^2 u_8}{2}\bigg) F^{(3)}[u_0]\\
+\bigg(\dfrac{u_2^2 u_3^2}{4} +\dfrac{u_1 u_3^3}{6} +\dfrac{u_2^3 u_4}{6} +u_1 u_2 u_3 u_4+\dfrac{u_1^2 u_4^2}{4}\\ +\dfrac{u_1 u_2^2 u_5}{2}+\dfrac{u_1^2 u_3 u_5}{2}+\dfrac{u_1^2 u_2 u_6}{2} +\dfrac{u_1^3 u_7}{6}\bigg)F^{(4)}[u_0] \\
+\bigg(\dfrac{u_2^5}{120}+\dfrac{u_1 u_2^3 u_3}{6} +\dfrac{u_1^2 u_2 u_3^2}{4} +\dfrac{u_1^2 u_2^2 u_4}{4} +\dfrac{u_1^3 u_3 u_4}{6}+\dfrac{u_1^3 u_2 u_5}{6}+\dfrac{u_1^4 u_6}{24}\bigg) F^{(5)}[u_0]\\
+\bigg(\dfrac{u_1^2 u_2^4}{48} +\dfrac{u_1^3 u_2^2 u_3}{12} +\dfrac{u_1^4 u_3^2}{48} +\dfrac{u_1^4 u_2 u_4}{24} +\dfrac{u_1^5 u_5}{120}\bigg) F^{(6)}[u_0]\\
+\bigg(\dfrac{u_1^4 u_2^3}{144} +\dfrac{u_1^5 u_2 u_3}{120} +\dfrac{u_1^6 u_4}{720}\bigg) F^{(7)}[u_0]+
\bigg(\dfrac{u_1^6 u_2^2}{1440}+\dfrac{u_1^7 u_3}{5040}\bigg) F^{(8)}[u_0]\\
+\dfrac{u_1^8 u_2 }{40320}F^{(9)}[u_0]+\dfrac{u_1^{10} }{3628800}F^{(10)}[u_0]\,.\)

\bibliographystyle{cmes}
\bibliography{azreg-5-9}

\end{document}